\documentclass[article]{IEEEtran}
\IEEEoverridecommandlockouts
\usepackage{amsmath,amssymb,amsfonts}
\usepackage{algorithmic}
\usepackage{graphicx}
\usepackage{textcomp}
\usepackage{xcolor}
\def\BibTeX{{\rm B\kern-.05em{\sc i\kern-.025em b}\kern-.08em
    T\kern-.1667em\lower.7ex\hbox{E}\kern-.125emX}}

\usepackage[T1]{fontenc}
\usepackage[utf8]{inputenc}

\usepackage{url, enumerate}
\usepackage{hyperref} 
\usepackage{verbatim}
\usepackage[francais,english]{babel}
\newcommand{\qed}{\bigstar} 

\newcommand{\base}{{\rm base}}
\newcommand{\cf}{{\rm cf}}

\newcommand{\supp}{{\rm supp}}

\newcommand{\ltr}{\!\triangleleft}

\newcommand{\ii}{\mathit i}
\newcommand{\Pset}{\mathcal P}

\newenvironment{proof}{\noindent{\bf Proof.}}{\par\bigskip}

\newenvironment{proof-}{\noindent{\bf Proof}}{\par\bigskip}

\newtheorem{THEOREM}{Theorem}[section]

\newtheorem{Conclusion}[THEOREM]{Conclusion}

\newtheorem{Hypothesis}[THEOREM]{Hypothesis}

\newtheorem{LEMMA}[THEOREM]{Lemma}

\newtheorem{Main Theorem}[THEOREM]{Main Theorem}
\newenvironment{main Theorem}{\begin{Main Theorem}} 
{\end{Main Theorem}}

\newtheorem{Theorem}[THEOREM]{Theorem}
\newenvironment{theorem}{\begin{Theorem}}{\end{Theorem}}

\newtheorem{Definition}[THEOREM]{Definition}
\newenvironment{definition}{\begin{Definition}}{\end{Definition}}

\newtheorem{Conventions}[THEOREM]{Conventions}

\newtheorem{Main Definition}[THEOREM]{Main Definition}
\newenvironment{main definition}{\begin{Main Definition}}
{\end{Main Definition}}

\newtheorem{Lemma}[THEOREM]{Lemma}
\newenvironment{lemma}{\begin{Lemma}}{\end{Lemma}}

\newtheorem{Notation}[THEOREM]{Notation}

\newtheorem{Convention}[THEOREM]{Convention}

\newtheorem{Observation}[THEOREM]{Observation}
\newenvironment{observation}{\begin{Observation}}
{\end{Observation}}

\newtheorem{Main Fact}[THEOREM]{Main Fact}
\newenvironment{main Fact}{\begin{Main Fact}}{\end{Main Fact}}

\newtheorem{Fact}[THEOREM]{Fact}

\newtheorem{Subfact}[THEOREM]{Subfact}

\newtheorem{Claim}[THEOREM]{Claim}
\newenvironment{claim}{\begin{Claim}}{\end{Claim}}

\newtheorem{Main Claim}[THEOREM]{Main Claim}
\newenvironment{main claim}{\begin{Main Claim}}{\end{Main Claim}}

\newtheorem{Crucial Claim}[THEOREM]{Crucial Claim}
\newenvironment{crucial claim}{\begin{Crucial Claim}}{\end{Crucial Claim}}

\newtheorem{Subclaim}[THEOREM]{Subclaim}

\newtheorem{Sublemma}[THEOREM]{Sublemma}

\newtheorem{Corollary}[THEOREM]{Corollary}
\newenvironment{corollary}{\begin{Corollary}}{\end{Corollary}}

\newtheorem{Example}[THEOREM]{Example}
\newenvironment{example}{\begin{Example}}{\end{Example}}

\newtheorem{Conjecture}[THEOREM]{Conjecture}
\newenvironment{conjecture}{\begin{Conjecture}}{\end{Conjecture}}

\newtheorem{Discussion}[THEOREM]{Discussion}

\newenvironment{Proof of the Subfact}
{\noindent{\bf Proof of the Subfact.}}{\par\bigskip}

\newenvironment{Proof of the Theorem}
{\noindent{\bf Proof of the Theorem.}}{\par\bigskip}

\newenvironment{Proof of the Proposition}
{\noindent{\bf Proof of the Proposition.}}{\par\bigskip}

\newenvironment{Proof of the Conclusion}
{\noindent{\bf Proof of the Conclusion.}}{\par\bigskip}

\newenvironment{Proof of the Observation}
{\noindent{\bf Proof of the Observation.}}{\par\bigskip}

\newenvironment{Proof of the Fact}
{\noindent{\bf Proof of the Fact.}}{\par\bigskip}

\newenvironment{Proof of the Lemma}
{\noindent{\bf Proof of the Lemma.}}{\par\bigskip}

\newenvironment{Proof of the Claim}
{\noindent{\bf Proof of the Claim.}}{\par\bigskip}

\newenvironment{Proof of the Corollary}
{\noindent{\bf Proof of the Corollary.}}{\par\bigskip}

\newenvironment{Proof of the Subclaim}
{\noindent{\bf Proof of the Subclaim.}}{\par\medskip}

\newenvironment{Proof of the Main Claim}
{\noindent{\bf Proof of the Main Claim.}}{\par\bigskip}

\newenvironment{Proof of the Crucial Claim}
{\noindent{\bf Proof of the Crucial Claim.}}{\par\bigskip}

\newcommand{\into}{\rightarrow}

\newcommand{\rest}{\upharpoonright}  

\newcommand{\deq}{\buildrel{\rm def}\over =}

\newcommand{\DS}{D\!S}
\newcommand{\D}{\Delta}

\newcount\skewfactor
\def\mathunderaccent#1#2 {\let\theaccent#1\skewfactor#2
\mathpalette\putaccentunder}
\def\putaccentunder#1#2{\oalign{$#1#2$\crcr\hidewidth
\vbox to.2ex{\hbox{$#1\skew\skewfactor\theaccent{}$}\vss}\hidewidth}}

\newcommand{\dom}{\mbox{\rm dom}}

\begin{document}

\title{On the ABK Conjecture, alpha-well Quasi Orders and Dress-Schiffels Product\\
\thanks{Mirna D\v zamonja's research was supported through the ERC H2020 project FINTOINF-MSCA-IF-2020-101023298. She thanks the Milyon laboratory for their hospitality in Lyon in December 2019 and the centre IHSPT at the Université Panthéon-Sorbonne, Paris, where she is an Associate Member. We are grateful to the referee for his detailed and helpful remarks, including finding inaccuracies to correct.}
}

\IEEEaftertitletext{\fbox{In loving memory of our co-author Maurice Pouzet, who died on the 31st of December, 2023.}}

\author{\IEEEauthorblockN{
Uri Abraham}
\IEEEauthorblockA{\textit{Department of Mathematics} \\
\textit{Ben-Gurion University}\\
84105 Beer-Sheva, Israel \\
abraham@math.bgu.ac.il}

\IEEEauthorblockN{
Robert Bonnet}
\IEEEauthorblockA{\textit{Laboratoire de math\'{e}matiques}\\
\textit{Universit\'{e} de Mont Blanc-Savoie}\\
73376 Le Bourget-du-Lac cedex, France\\
rbonnet.perso@gmail.com}

\IEEEauthorblockN{
Mirna D\v zamonja}
\IEEEauthorblockA{\textit{IRIF (CNRS)} \\
\textit{Université de Paris Cité}\\
75205 Paris cedex 13,
France \\
mdzamonja@irif.fr}

\IEEEauthorblockN{
\fbox{Maurice Pouzet}}
\IEEEauthorblockA{\textit{D{\'e}partement de math\'{e}matiques} \\
\textit{Universit\'{e} Claude Bernard, Lyon 1}\\ 69622 Villeurbanne cedex,  France\\
maurice.pouzet@univ-lyon1.fr}}

\maketitle

\begin{abstract}
The following is a 2008 conjecture from \cite{bqoconjecture}\footnote{We named the conjecture by the initials of the authors of \cite{bqoconjecture}.}

[ABK Conjecture] Every well quasi order (wqo) is a countable union of better quasi orders (bqo).

We obtain some partial progress on the conjecture, in that
we show that the class of orders that are a countable union of better quasi orders (sigma-bqo) is closed under various operations. These include diverse products, such as the little known but natural Dress-Shieffels product. We develop various properties of the latter. In relation with the main question, we explore the class of
alpha-wqo for countable ordinals alpha and obtain several closure properties and a Hausdorff-style classification theorem. Our main contribution is the discovery of various properties of sigma-bqos and ruling out potential counterexamples to the ABK Conjecture. 
\end{abstract}


\section{Introduction}\label{introduction} The notion of a well-quasi order is ubiquitous in mathematics, by the naturalness of its definition and by the many applications
it has found in all areas of mathematics. The early history is well described in Joseph Kruskal's article \cite{Kruskal} and later references abound. The notion has had
another life in theoretical computer science. For example, 
in 2017 the influential Computer Aided Verification (CAV) award was given to Parosh Abdulla, Alain Finkel, Bengt Johnsson and Philippe Schnoebelen, who showed that the notion of {\em well-quasi ordering (wqo)} can be used to single out a large class of infinite-state transition systems that have significant decidability properties. However, the class of wqo orders lacks some closure properties, such as the closure under infinite products. This has led to the definition of stronger notions which do have closure properties, notably the notion of {\em better quasi ordering (bqo)}. This paper is about the relation of the two notions, which we shall now introduce.

We define the notion of wqo, passing through a translation to the more familiar notion of partially ordered sets or posets.
Since every quasi-order can be considered as a poset after a simple factoring by $\le \cap \ge$, in this article we simply work with wqo posets usually denoted by $(P,\le_P)$. This does not lead to any loss of generality in whatever we have to say. 

A poset is a {\em wqo} if it has no infinite antichains or strictly decreasing sequences. It is of interest to know which natural operations with posets preserve the property of being a wqo.

A {\em restriction} of a poset $P$ is a subset of $P$ with the induced suborder. It follows that the notion of wqo is preserved by restrictions. 
Finite products of posets are defined by considering coordinatwise order on the Cartesian product of the posets in question. It is well known and easily proved that a finite product of wqo orders is a wqo.

For infinite products of $\alpha$ factors we use a more refined notion, which in fact reduces to the usual product in the case of $\alpha$ finite. Throughout the paper
$\alpha$ stands for an ordinal.

\begin{definition}\label{def:prod} Suppose that $\alpha$ is an ordinal and $P$ a poset. By $P^\alpha$ we denote the factor poset obtained from the quasi order on the set ${}^\alpha P$
of $\alpha$-sequences from $P$ with the {\em embeddability order}, given by:
\[
s\le t
\mbox{ iff there is a strictly increasing }\rho:\,\alpha\to\alpha\]
\[\mbox{ such that }s(\gamma)\le_P t(\rho(\gamma)), \mbox{ for all }\gamma<\alpha.
\]
(Hence the factoring is done through the equivalence relation $s\sim t$ iff $s\le t$ and $t\le s$.) Similarly, by $P^{<\alpha}$ we denote the factor poset of the set of all sequences of length $<\alpha$ from $P$, ordered by embeddability.  
\end{definition}

The infinite product $P^\omega$ for a wqo order $P$ is not always a wqo, as was proved by Richard Rado in \cite{rado54} who gave a counterexample now known as the Rado order $R$. Crispin Nash-Williams in \cite{bqo}
defined the notion of a {\em better quasi order (bqo)}, which a strengthening of the notion of wqo and is specifically designed so not to be satisfied by $R$. Nash-Williams proved that for a bqo $P$, the product $P^\omega$ is a wqo, where $\omega$ is the first infinite ordinal. A full characterisation of bqo $P$ in terms of the powers 
of $P$ being wqo was given by 
Maurice Pouzet \cite{Pozetchar}, improving on the results of Nash-Williams and stating:

\begin{theorem}[Nash-Williams/Pouzet] \label{Mauricecharacterisation1}A wqo $P$ is a bqo iff $P^{<\omega_1}$ is wqo.
\end{theorem}

The above characterisation naturally leads to the study of wqos that satisfy that the power $P^\alpha$ is wqo, for an ordinal $\alpha\in[\omega, \omega_1)$. It turns out that for indecomposable $\alpha$, a precise characterisation can be given in terms of the Ramsey-theoretic notion of $\alpha$-wqo, whose definition we recall in \S \ref{reviewbqo}. Then the following theorem follows from the work of Pouzet \cite{Pozetchar} and of Alberto Marcone in \cite{marcone94}.

\begin{theorem}[Marcone/Pouzet]\label{alphawqo} Let $\alpha <\omega_1$ be an indecomposable ordinal. A poset $P$ is {\em $\alpha$-wqo} iff the poset $P^{<\alpha}$ is wqo.
\end{theorem}

As mentioned in the above, the notion of bqo was developed by analysing the properties of Rado's example and finding a property that this example does not satisfy. Hence Rado's order is an example of a wqo order which is not bqo. Further, uncountable examples exist, often obtained by generalising Rado's example to barriers (see below)
and differentiating $\alpha$-wqos from $\beta$-wqos, for $\beta>\alpha$.
It does seem that there are not that many natural examples of such orders, in the sense that the wqo theory is often applied to show that some particular order does not have an infinite decreasing sequence (see for example the celebrated solution by Richard Laver of the Fraïssé conjecture \cite{laver71}), but no such applications have so far given a wqo which is known not be bqo. Nash-Williams states on 
p.700 of \cite{zbMATH03232673} that ``one  is inclined to  conjecture  that  most  w.q.o. sets which  arise in a reasonably  `natural'  manner  are likely to be b.q.o'' . 
A possible combinatorial link between the notions of wqo and bqo is proposed by
the following conjecture of Uri Abraham, Robert Bonnet and Wies\l aw Kubi\'s in \cite{bqoconjecture}. (It was named the ABK Conjecture by D{\v z}amonja and Pouzet.)

\begin{conjecture}[ABK Conjecture]\label{AB} Every wqo is a countable union of bqo.
\end{conjecture}

By a countable union, we mean a union $\bigcup_{n<\omega}P_n$ in which each $P_n$ is a bqo. No requirement on the relations between the elements of
different $P_n$s is assumed, so we can equivalently consider disjoint or increasing unions - the only requirement is that there is some splitting of the underlying
set into countably many pieces so that the restriction of the order on each peace is a bqo. In particular, every countable order is a countable union of bqos. We 
are only interested in the orders that themselves are wqo. Hence
we shall define {\em $\sigma$-bqo}s as those wqos which are countable unions of bqos. That is, the ABK Conjecture states that every wqo is $\sigma$-bqo. 

ABK Conjecture rules out the known examples of wqos which are not bqos. For example, Rado's order is $\sigma$-bqo because it is countable. 

If the conjecture is true, it would allow us to study closure properties of wqos which have not previously been noticed. Namely, wqos are notoriously not preserved when it comes to various operations, such as infinite products, while the bqos are. In this paper we show that the property of being the countable union of bqos is also preserved by many operations. Of course, even if ABK Conjecture were to be true, this would not help us with the closure properties of wqos, since they are known not to hold. However, it would give us a handle of analysing and classifying the failure of these properties. For example, if the conjecture were to be true, even when we know that we have a wqo whose infinite product is not a wqo, we would still be able to conclude that we can represent the product as a countable union of nice subsets, each of which is a bqo. As infinite products are typically uncountable, this would represent progress. Infinite products are studied in the combinatorics of words and in fact this was the motivation of Rado and Nash-Williams, to start with. The study of infinite products of wqos has also led to
an alternative to bqos, which are the Noetherian spaces studied by Jean Goubault-Larrecq in \cite{Noetherian} and elsewhere.

A poset $P$ is countably indivisible 
if whenever it is partitioned into countably many pieces, then $P$ is isomorphic to a subposet of one of these pieces. This is a natural generalisation of the notion of
indivisibility. The latter notion has been studied in various contexts, including relational structures (see \cite{zbMATH06256128} for definitions, further references and historical remarks) as well as ultrametric topological spaces (\cite{zbMATH05306960}). There does not seem to be much knowledge available about countably indivisible structures.
An interesting application of the positive answer to the ABK Conjecture would be the conclusion
that every countably indivisible wqo is in fact a bqo. 

These have been our motivations to study not just the conjecture, but also the property of being $\sigma$-bqo in general. 
Our contribution has been to rule out those counterexamples to the ABK Conjecture that one could imagine constructing by any operation on wqos that appears in the literature, at least to the extent of our knowledge. This does not mean that a counterexample does not exist, but it does seem to indicate that if it is to be found a really new idea would be needed. 

Our research of the literature led us to find what we consider a very interesting work, relating to what we have named Dress-Schiffels product.
This special kind of product of posets was defined by Andreas Dress and Gerhard Schiffels in
a technical report \cite{DressSchiffels}. A variant of the product appears in \S8 of \cite{MR1440456}. Apart from these hard-to-find references that appeared
30 or more years ago, the product does not exist in any published literature. Yet, it seems like a very natural way to define the finite support product of posets (to handle the fact that the orders are not linear, the definition needs an extra twist). In this paper we revived this way of taking products, which we feel is likely to find applications in the logic community, and we developed various properties of such products. In particular, we proved that such a product preserves the properties of being $\alpha$-wqo and $\sigma$-bqo (and so in particular cannot produce counterexamples to the ABK Conjecture).

As a final note, let us mention that in the last stages of preparation of this article we discovered a survey paper \cite{Mauricesurvey3} that our co-author Maurice Pouzet published posthumously. The article not only shows the depth of the general area of research of M. Pouzet, but also offers a precious source of interesting open questions and directions. Hence it is our pleasure to recommend it to any reader interested in wqos and all other relational structures.

\section{A Brief Presentation of our Results}\label{sec:present_results}
In this paper we show the partial progress that we have obtained, by eliminating various natural candidates for a counterexample. On the way, we have developed quite a lot of theory of $\alpha$-wqo and the $\sigma$-bqo orders, which we believe is of independent interest.

Our main results, enumerated as in the main text of the paper are:

\smallskip 

{\bf Theorem \ref{alpha<alpha}}. {\em Suppose that $P^\alpha$ is $\varepsilon$-wqo for some
$\varepsilon\in [\omega, \omega_1)$. Then so is $P^{<\alpha}$}.

\smallskip 

{\bf Theorem \ref{newknowledge}}.{\em Suppose that $\alpha>0$ is a limit ordinal. If
 $I(P^{<\alpha})$ is wqo, then so are $P^\alpha$ and $P^{<\alpha}$.}
 
 \smallskip 

 Here, for the order $Q$, the notation $I(Q)$ stands for the poset consisting of the initial segments of $Q$ ordered by inclusion.

{\bf Theorem \ref{products}}. {\em Let $1\le\alpha < \omega_1$. Suppose that $P$ is $\sigma$-bqo and that 
$P^\alpha$ is wqo. Then the orders $P^{<\alpha}$, $I(P^{<\alpha})$ and $P^\alpha$ are all $\sigma$-bqo.}

\smallskip 
 
{\bf Theorem \ref{lexsums}} Suppose that $Q$ and $\langle P_q:\,q\in Q\rangle$ are $\sigma$-bqo orders.
Then so is the lexicographic sum 
\[
\Sigma_{q\in Q} P_q.
\]

\smallskip 

In \S \ref{sec:DS} we recall the Dress-Schiffels product $\otimes^{\DS}_{i\in I} P_i$ and prove

\smallskip 

{\bf Theorem \ref {generalisedDressShiffels}}. If $I$ and $P_i$ for every $i\in I$ are bqo ($\alpha$-wqo for some $\alpha\in [\omega,\omega_1)$, $\sigma$-bqo) posets, then so is $\otimes^{\DS}_{i\in I} P_i$.

\smallskip 

A poset $P$ is said to be {\em FAC} if all its antichains are finite. This is therefore a weakening
of the notion of wqo, which is nevertheless interesting in applications, see notably \cite{lmcs:3928} where it is shown that this requirement on the underlying order of a transition system, implies the decidability of the covering problem. In \S\ref{sec:Hausdorff} we introduce a naturally defined notion of {\em $\alpha$-FAC} orders
and show that the class of such orders is obtained by closing the class of $\alpha$-wqos under some naturally defined operations on orders. This is a Hausdorff-style classification theorem for {\em $\alpha$-FAC} orders. Although it is not directly connected to the ABK Conjecture, we have felt that the theorem fits well in a paper that develops the theory of $\alpha$-wqos. We need a bit of notation in order to explain this theorem.

\smallskip 

{\bf Definition \ref{reasonable}}. A class $\mathcal{G}_0$ of linear orderings is {\em reasonable}
    if and only if $\mathcal{G}_0$ contains a nonempty ordering and
   is closed under reversals and restrictions.
   
\smallskip 

{\bf Definition \ref{reasonable2}}. Given a reasonable class $\mathcal{G}_0$ of linear orderings,
    the $\alpha$-{\em closure} $\alpha cl( {\mathcal G}_0 )$ of $\mathcal{G}_0$ is the least
   class of posets which contains $\mathcal{G}_0$ and is closed
   under the operations:
\begin{itemize}
\item Lexicographic sum with index set  either an  $\alpha$-wqo poset, the inverse
   of a  $\alpha$-wqo poset, or an element of $\mathcal{G}_0$.  
\item Augmentation (having the same domain, but possibly adding relations).
\end{itemize}
 
\smallskip 

{\bf Theorem \ref{forbiddencor}}. Let $\mathcal{G}_0$ be a reasonable class of posets and
 let $\mathcal{G}$ be the class of $\alpha$-FAC posets such that every chain
  is in $\mathcal{G}_0$,
    and assume in addition that:
\begin{enumerate}
\item $\mathcal{G}_0$ contains all well-orderings, and is closed under
   lexicographic sums with index set in $\mathcal{G}_0$.
\item $\mathcal{G}$ is closed under augmentations.
\end{enumerate}

 Then $\mathcal{G} = \alpha cl(\mathcal{G}_0)$.
 
 \section{Organisation of the Paper}\label{sec:organisation}
\S \ref{reviewbqo} contains a review of the Nash-Williams theory of barriers and better quasi order. The techniques and the results of this section are used in many places in our proofs. This theory is well known, but not considered an easy read. We hope that the way that we present it will change that perspective and convince the reader of the usefulness of the concept. In this section we also introduce one of our main topics, the notion of $\alpha$-wqo. The short \S \ref{techniques} introduces some easily observed but very useful techniques which are isolated here so as to be referred to in the main body of the paper. Our results start in \S\ref{sec: connections}, where we present some of the deep connections that exist between the level of well quasi ordering on the products of the type $P^\alpha$ and $P^{<\alpha}$, as well as the poset of initial segments $I(P^{<\alpha})$. \S 
\ref{sec:products} considers the $\alpha$-wqo and $\sigma$-bqo properties of products of type $P^\beta$, also using some of the results of 
\S\ref{sec: connections}. The proofs of the theorems in these two sections are some of the main proofs in the paper. Section \ref{sec:lexsums} is quite light, it discusses the closure of the various notions in question under lexicographic sums and the proofs are readily obtained. Another main section is \S\ref{sec:DS},
where we consider a special product of orders. It is the Dress-Schiffels product and one of our main theorems is that the notions of $\alpha$-wqo, bqo and $\sigma$-bqo are closed under this product. Dress-Schiffels product in itself deserves more attention, as almost nothing about it is available in the published literature. The final section 
\S\ref{sec:Hausdorff} contains a Hausdorff-style build up theorem which shows how to obtain the class of $\alpha$-FAC orders by building up from the class of $\alpha$-wqo orders. The proofs are of a certain length, however built on techniques explored in \cite{abcdzt}. 

The working part of the paper ends with Section \ref{sec:proofs} which contains the proofs that we considered better presented after the main body of the paper.
This section is followed by a Conclusion in \S\ref{sec: conclusion}.

\section{Better Quasi Orders and the Definition of $\alpha$-wqo}\label{reviewbqo} In this section we recall the main points of Nash-Williams' theory of better quasi orders, bqo.
This theory dates from 1965 \cite{NashWilliamsoriginal} and it has proved its usefulness in many contexts. Yet it is still often perceived as too technical and lacking intuition. We give a short presentation which shows how this theory naturally emerges from Ramsey's theorem and how in the end, the notion of a bqo is indeed a natural extension of that of a wqo. We mention the advantages of bqos, for example the closure of the class of bqos under various operations. Finally, we open the way for the introduction of the notion of $\alpha$-wqo, which is done at the end of the section.

We say that a sequence $\langle p_i:\,i<\omega\rangle$ of elements of a poset $P$ is a {\em bad sequence} if there are no $i<j$ such that
$p_i\le p_j$.  By Ramsey's theorem,
a poset $P$ is wqo iff it does not allow for a bad sequence. A generalisation of this characterisation to the notion of barriers, Definition \ref{def: barrier}, is the combinatorial 
heart of Nash-Williams' development of bqo. Barriers are more complicated structures than the bare set of natural numbers,  but they turn out to have enough Ramsey-theoretic properties so to be useful in generalisation of the notion of a bad sequence.

We denote by $ [\omega]^{<\omega}$ the set of finite subsets of $\omega$.

\begin{definition}[Barrier]\label{def: barrier}
(1) For a set $A\subseteq [\omega]^{<\omega}$ we define the {\em base of} $A$ as $\base(A)=\bigcup A$.

{\noindent (2)} An infinite set $B\subseteq [\omega]^{<\omega}$ is a {\em barrier} if:
\begin{enumerate}
\item $B$ is an antichain under $\subseteq$, that is, for $x\neq y\in B$ it is not the case that $x\subseteq y$.
\item
The base $U$ of $B$ has the property that every infinite subset of $U$ has an initial  segment in $B$.
\end{enumerate}
{\noindent (3)} A {\em sub-barrier of $B$} is a barrier $C$ such that $C\subseteq B$. 
\end{definition}

The definition of a barrier is admittedly not intuitive. However, barriers are the heart of the bqo theory, as they allow for a beautiful Ramsey-type theory. We recall the main points of it, starting with the Nash-Williams' partition theorem. It is Theorem 1 from \cite{NashWilliamsoriginal}.\footnote{We mention that the original Nash-Williams partition theorem has been much extended through the work in descriptive theory, starting from Fred Galvin and Karel Prikry's work in \cite{GalvinPrikry}, where Nash-Williams becomes a corollary. However, we are interested in the purely combinatorial content of the bqo theory here, since this is what is used in the work on $\alpha$-wqo.}

\begin{theorem}[Nash-Williams partition theorem]\label{thm:Nash-Williams} If $B$ is a barrier and $B=B_0\cup B_1$ is a partition of $B$, then there is 
$i<2$ such that $B_i$ contains a barrier. 
\end{theorem} 

Obviously, by applying the theorem finitely many times, we obtain that in any finite partition $B_0, B_1,\ldots B_k$ of a barrier, there is a part $B_i$ which contains a barrier. In other words there exists a barrier $C\subseteq B_i$ whose base
$Y$ is a subset of $\bigcup B_i$. This already implies that $C=B\cap [Y]^{<\omega}$ and hence for every index $j\neq i$ 
we have $B_j\cap [Y]^{<\omega}=\emptyset$.

Barriers are particular subsets of  $[\omega]^{<\omega}$ and hence can be lexicographically ordered. A theorem of Pouzet
in \cite{Pozetchar}
(see Lemma 1.5 in \cite{marcone94} for a proof in English) gives that barriers are well ordered by the lexicographic order.
The order type of a barrier under the lexicographic order is called the {\em type} of the barrier. This notion will intervene in
Definition \ref{nashwilliams} and in \S\ref{sec:alpahwqo}, where we discuss
the $\alpha$-wqos.

\begin{example}\label{getRamsey} There are no barriers of finite order type. The simplest barrier is $\{\{n\}:\,n<\omega\}$. It is the only barrier of order type 
$\omega$. 
\end{example}

Nash-Williams partition theorem can be extended to a full Ramsey-style Theorem \ref{Ramseybar}, but we first need to introduce the powers of a barrier. To do so, we need the notion of shift-extension.

\begin{definition}\label{shiftextends} (1) The relation $\ltr$ of {\em shift-extension} in $ [\omega]^{<\omega}$ is defined as 
follows.
For $r,s\in[\omega]^{<\omega}\setminus\{\emptyset\}$, we let $r\ltr s$ if:
\begin{enumerate}
\item $\min(r) < \min(s)$ and
\item $r\setminus \{\min( r)\}$ is a proper
initial segment of $s$.
\end{enumerate}

{\noindent (2)} Let $U$ be a barrier. Define
\[ U^2= \{ s\cup t\mid s,t\in U \wedge s\ltr t \}.\]

\end{definition}

For example $\{4\} \ltr \{ 7,34, 45\}$, and $\{4,7,34\}\ltr \{ 7,34, 45\}$, but $\{4,7\}\not\ltr \{7\}$. 

Note that if $U$ is a barrier then $U^2$ is a barrier and its base is the same as the base of $U$, namely $\bigcup U^2=\bigcup U$. 
One should not confuse the barrier $U^2$ with the cartesian product $U\times U$. If $u\in U^2$ then there exist uniquely determined elements of $U$, $s$ and $t$, for which $u=s\cup t$. We write $u=[s,t]$.

An important lemma is:

\begin{lemma}\label{squarebarrier}
\label{LemSB} (1) If $U$ is a barrier, then it contains a strictly $\ltr$-increasing infinite sequence.

{\noindent (2)} Any sub-barrier of $U^2$ has the form $V^2$ for some sub-barrier $V$ of $U$.
\end{lemma}

For a barrier $U$ we can also define $U^3=\{u\cup v\cup w\mid u,v,w\in U \wedge u\,\ltr v\, \ltr w \}$.
 If $x,y\in U^2$ and $x\ltr y$, then $x=[a,b]$ and $y=[b,c]$ for some $a,b,c\in U$. It follows that $U^3 = (U^2)^2$.
Using the previous lemma we get that a sub-barrier of $U^3$ is of the form $V^3$ where $V$ is a sub-barrier
of $U$, and so on. 

These observations, all due to Nash-Williams in \S3 of \cite{zbMATH03232673} (although they are easier to follow in the survey paper \cite{Milner} by Eric Milner) allow us to formulate a version of Ramsey's theorem for barriers. It implies the classical Ramsey theorem when applied to the barrier consisting of singletons.

\begin{theorem}[Ramsey theorem for barriers]\label{Ramseybar} Suppose that $n,k\ge 1$ are natural numbers, $B$ is a barrier, and $c:\,B^n\to k$
is a function (often called a colouring). Then there exist a sub-barrier $H$ of $B$ and a colour $l^\ast<k$ such that 
\[
c\rest H^n:\,H\into \{l^\ast\},
\]
that is, $H^n$ is monochromatic.
\end{theorem}

We can finally make the connection with the notion of a bad sequence, through the idea of barrier sequences, or arrays.

\begin{definition}[Good, Bad and Perfect]\label{goodandbad} Let $U$ be a barrier, $P$ a poset and $f: U\to P$ a function (such a function is often called a $U$ {\em barrier
sequence} on $P$ and in some literature, an array). 

{\noindent (1)} We say that $f$ is {\em good} if for some $s,t\in U$ such that $s\ltr t$, $f(s)\leq_P f(t)$.
If $f$ is not good we say that it is {\em bad}. 

{\noindent (2)} We say that $f$ is {\em perfect} if for all $s,t\in U$ for which $s\ltr t$
we have that $f(s)\leq_P f(t)$.
\end{definition}

We can now state the following definition:

\begin{definition}[$\alpha$-wqo and bqo]\label{nashwilliams} Let $\alpha\in[\omega, \omega_1)$ and let $P$ be a poset. We say  that $P$ is {\em $\alpha$-wqo} iff for every
barrier $B$ with order type $\le\alpha$, every  map $f:\,B\to P$ is good with respect to $P$. The poset $P$ is {\em bqo} iff it is 
$\alpha$-wqo for all $\alpha<\omega_1$.
\end{definition}

In particular, Definition \ref{nashwilliams} gives that a poset is wqo iff it is $\omega$-wqo.

An important Ramsey-style theorem is given in the following. It is Lemma 19 in \cite{zbMATH03232673}
and its proof follows easily from what we stated already.

\begin{theorem}[Nash-Williams]\label{NW} Let $f$ be a function from a barrier $U$ into a poset $P$. Then there is a sub-barrier
$V\subseteq U$ such that the restriction of $f$ to $V$, $f\restriction V$ is either bad or perfect.
\end{theorem}

\begin{proof} Let $E=\bigcup U$ be the base of $U$ (and of $U^2$). Partition $U^2$ into two sets by the function $g:U^2\to 2$ 
defined as $g(s\cup t)=0$ iff $f(s)\leq_P f(t)$ (for
$[s, t]\in U^2$). By the Ramsey theorem for barriers, there is a sub-barrier $V$ of $U$ such that $g$ is constant on $V^2$.
If $g$ has the constant value $0$ on $V^2$ then $f\restriction V$ is a perfect barrier sequence, and otherwise it is a bad barrier 
sequence.
 $\qed_{\ref{NW}}$
\end{proof}

A consequence of Nash-Williams' partition results for the barriers is that the class of bqo orders is closed under 
various operations. For example (for the proof see (2.15), pg. 496 of \cite{Milner}):

\begin{theorem}\label{bqoproducts} If $P$ and $Q$ are bqo, so is the product $P\times Q$.
\end{theorem}

In the main body of the paper we shall discuss various other closure properties of both bqo and $\alpha$-wqo orders.

 \subsection*{A few more words on $\alpha$-wqo}\label{sec:alpahwqo} Possible order types of barriers are a proper subset of the set of countable ordinals, as follows from the following theorem combining several results of Marc Assous in \cite{MR366758}:
 
 \begin{theorem}[Assous]\label{Pozetchar} The order type of any barrier is of the form $\omega^\beta\cdot n$, where $\beta>0$ is a countable ordinal and $n>0$ a natural number. If $\beta<\omega$, then $n=1$. Moreover, for each ordinal of the above type, there is a barrier whose order type is that ordinal.
\end{theorem}

Assous' result explains that the most interesting type of ordinals $\alpha$ for the notion of $\alpha$-wqo are the indecomposable $\alpha$.
For such $\alpha$ there is a characterisation, which follows from the work of Pouzet \cite{Pozetchar} and that of Alberto Marcone in \cite{marcone94}:

\begin{theorem}[Marcone/Pouzet]\label{alphawqo1}Let $\alpha <\omega_1$ be an indecomposable ordinal. A poset $P$ is {\em $\alpha$-wqo} iff the poset $P^{<\alpha}$ is wqo.
\end{theorem}

Since the only barrier of order type $\omega$ is the one consisting of singletons, we easily obtain that a poset $P$ is $\omega$-wqo iff it is wqo, hence $\omega^2$ is the first ordinal $\alpha$ such that the class of $\alpha$-wqos is not immediately reduced to the class of wqos.

\section{Some Useful Techniques and Remarks}\label{techniques}
Several proofs in the paper will rely on applications of certain techniques, which we detail here.

Definition \ref{nashwilliams} has several immediate consequences.

\begin{observation}\label{fromalphatobeta} (1) Suppose that $\beta\le\alpha<\omega_1$ and $P$ is $\alpha$-wqo. Then $P$ is
$\beta$-wqo.

{\noindent (2)} Suppose that $P$ is bqo ($\alpha$-wqo) and $Q$ is a suborder of $P$. Then $Q$ is bqo ($\alpha$-wqo).

{\noindent (3)} Suppose that $P$ is $\sigma$-bqo and $Q$ is a suborder of $P$. Then $Q$ is $\sigma$-bqo.
\end{observation}

\begin{proof} (1) Immediate, since a barrier of order type $\le\beta$ is also of order type $\le\alpha$.

{\noindent (2)} Suppose that $U$ is a barrier (of order type $\le\alpha$) and $f:\,U\to Q$ is given. In particular, $f$
is a $U$ barrier sequence over $P$, so there are $s\ltr t\in U$ such that $f(s)\leq_P f(t)$ and hence $f(s)\leq_Q f(t)$.

{\noindent (3)} This follows from (2), since if $P=\bigcup_{n<\omega} P_n$ where each $P_n$ is a bqo, then
$Q=\bigcup_{n<\omega} (P_n\cap Q)$ and each $P_n\cap Q$ is bqo by (2).
$\qed_{\ref{fromalphatobeta}}$
\end{proof}

\begin{lemma}\label{orderpreserving} Suppose that $P$ and $Q$ are posets and $1\le\alpha<\omega_1$ an ordinal.

{\noindent (1)} If $P$ is $\alpha$-wqo (bqo, $\sigma$-bqo) and there is a surjection $f:\,P\to Q$ satisfying
\[
p\le_P p' \implies f(p)\le_Q f(p'),
\]
for all $p,p'\in P$, then $Q$ is $\alpha$-wqo (bqo, $\sigma$-bqo).

{\noindent (2)} If $Q$ is $\alpha$-wqo (bqo, $\sigma$-bqo) and there is $g:\,P\to Q$ satisfying
\[
g(p)\le_Q g(p') \implies p\le_P p', 
\]
for all $p,p'\in P$, then $P$ is  $\alpha$-wqo (bqo, $\sigma$-bqo).

[Note that since we are dealing with posets, rather than with quasi-orders, such a $g$ is necessary an injection].
\end{lemma}

\begin{proof} We prove only the instance relating to $\alpha$-wqo, since the proofs in the case bqo and $\sigma$-bqo are similar.

(1) Let $f:\,P\to Q$ be a surjection satisfying the hypothesis.
Suppose further that $P$ is $\alpha$-wqo and that $h:B\to Q$ is a function from some barrier $B$ of order type $\le\alpha$ to
$Q$. For each $q\in Q$ choose an element $p_q$ of the set $f^{-1}(q)$. Then $h':B\to P$ given by $h'(s)=p_q$ iff
$h(s)=q$ is well-defined and therefore there exist $s\ltr t$ in $B$ with $h'(s)\le_P h'(t)$. Hence $h(s)\le_Q h(t)$ by the
definition of $h'$ and the requirement on $f$, and so $h$ is good.

(2) Let $g:\,P\to Q$ be a function satisfying the hypothesis.

Suppose that  $Q$ is a $\alpha$-wqo and that $h:\,B\to P$ is a function from some barrier $B$ of order type $\le\alpha$ to
$P$. Then $g\circ h:\,B\to Q$ and therefore there are $s\,\ltr t$ in $B$ with 
\[(g\circ  h)(s)\le_Q (g\circ  h)(t).
\]
By the hypothesis on
$g$, it follows that $h(s)\le_P h(t)$ and so $h$ is good.
$\qed_{\ref{orderpreserving}}$
\end{proof}

A surjection such as in Lemma \ref{orderpreserving}(1) is said to be {\em order-preserving}. A function 
such as in Lemma \ref{orderpreserving}(2) will be called {\em order-generating}.

As a side remark, we note that the above technique does not necessarily apply to the case of well founded posets. Namely, the existence of an order-generating function from $P$ to a well-founded $Q$ does not necessarily 
imply that $P$ is well-founded. For example, $Q$ could be an infinite countable antichain and $P$ the order $\omega^\ast$. This example demonstrates the difference between the notions of wqo and being well founded.

Let us also mention that there is not very much that one can say about the ABK Conjecture by using the combinatorics of infinite cardinals. The following simple remark is a rare one to be made by arguing about cardinalities, but it is perhaps worth
noting here among the introductory results. At least to some extent it explains why the naive approach of building a counterexample cannot work.

\begin{Observation}\label{countablecof} Suppose that $P$ is a wqo satisfying $|P|=\kappa$, where $\cf(\kappa)=\aleph_0$ and such that every $Q\subseteq P$ with $|Q|<\kappa$ is $\sigma$-bqo. Then $P$ is $\sigma$-bqo.
\end{Observation}

\begin{proof} Clearly, if $\kappa=\aleph_0$, the order $P$ is $\sigma$-bqo since every singleton is a bqo. Suppose that 
$\kappa>\cf(\kappa)=\aleph_0$ and let $\langle \kappa_n:\,n<\omega\rangle$ be a sequence of cardinals $<\kappa$
with $\sup_{n<\omega} \kappa_n=\kappa$. 

We decompose $P$ into a disjoint union $\bigcup_{n<\omega} Q_n$ where $| Q_n|=\kappa_n$. Therefore each $Q_n$
can be decomposed into a countable union of bqo, say $Q_n=\bigcup_{i<\omega} P^n_i$. Hence
$P=\bigcup_{i,n<\omega} P^n_i$, is a countable union of bqo.
$\qed_{\ref{countablecof}}$
\end{proof}

\section{Connections Between $P^\alpha$ and $P^{<\alpha}$ and Between $P$ and $I(P)$}\label{sec: connections}
Section \ref{sec:products} will consider the status of the ABK Conjecture in the case of products of the
form $P^\alpha$ and $P^{<\alpha}$. In order to approach these, we develop several structural properties of such orders in general. Throughout this section, $\alpha$ and $\beta$ stand for non-zero countable ordinals and $P$ for a poset.

An {\em initial segment} of a poset $P$ is a subset closed under $\le_P$, that is a subset $A$ of $P$ such that for all $y\in A$, if $x\le_P y$, then $x\in A$. For a poset $P$ we denote by $I(P)$ the set of all initial segments of $P$ partially ordered by the subset relation.

The main objects of this section are posets of the form $P^\alpha$, $P^{<\alpha}$ and $I(P^{<\alpha})$.

\begin{theorem}\label{alpha<alpha} Suppose that $P^\alpha$ is $\varepsilon$-wqo for some
$\varepsilon\in [\omega, \omega_1)$. Then so is $P^{<\alpha}$.
\end{theorem}

\begin{proof} The case $\alpha=0$ is trivial, so assume that $\alpha>0$. Since $P^\alpha$ is wqo, Observation \ref{fromalphatobeta}(1) implies that $P$ is wqo. Therefore $P$ has a minimal element, say $a$. 

For $f\in P^{<\alpha}$ let $\gamma_f\leq \dom(f)$ be the least ordinal $\gamma$ such
that $\forall \beta \in [\gamma,\dom(f))\ (f(\beta)=a)$. Thus, $\gamma_f = \dom(f)$ iff
$\forall \gamma<\dom(f)\ \exists \zeta \ (\gamma\leq\zeta<\dom(f) 
\,\wedge\, f(\zeta)\neq a)$.  Therefore $f\restriction [0,\gamma_f)$ has an unbounded cofinal set whose $f$ values are $\neq a$. 

The order-type of the interval $[\gamma_f, \dom(f))$  is
denoted $\tau(f)$. So $0\leq \tau(f)\leq \dom(f)<\alpha$. For example, if $f$ is the constant $a$ function
then $\gamma_f=0$ and $\tau(f)=\dom(f)$, and the other extreme case is when $f$ has no non-empty final
segment of points that have value $a$, and in that case $\tau(f)=0$.

The interval $[0,\gamma_f)$ is said to be the left interval of $f$, and $[\gamma_f, \dom(f))$ is its right interval. So $\tau(f)$ is the order-type of the right interval of $f$.

Suppose $B$ is a barrier of order-type $\leq \varepsilon$ and let $h:B\to P^{<\alpha}$ be
a given barrier sequence on $P^{<\alpha}$. We shall find $s\ltr t$ in $B$ such that $h(s) \leq_{P^{<\alpha}} h(t)$.

Consider the function $\tau\circ h: B \to \alpha$, where $\tau$ is the function described above.
Using Ramsey's theorem for barriers (\ref{Ramseybar}), we get a sub-barrier $C$ of $B$ such that either
for all $s\ltr t$ that are both in $C$, $\tau(h(s)) > \tau(h(t))$ or else
for all such $s\ltr t$, $\tau(h(s)) \leq \tau(h(t))$. The first option is ruled out
since Lemma \ref{squarebarrier}(1) gives us an infinite $\ltr$-increasing sequence in $C$ and hence the first option would give us an infinite descending sequence of ordinals.
Therefore, we have a sub-barrier
$C$ over which the second option holds: for all $s\ltr t$ that are both in $C$,
\begin{equation}
\label{Eq:1}
\tau(h(s)) \leq \tau( h(t)).
\end{equation}  

We define a function $\varphi:\,P^{<\alpha}\to P^\alpha$ by letting $\varphi(f)= f^\frown
 a_{[\dom(f),\alpha)}$, where $a_{[\dom(f),\alpha)}$ denotes the constant $a$ valued function defined on the interval $[\dom(f),\alpha)$.
Then we define $\varphi\circ h : C \to P^\alpha$, a barrier sequence over $P^\alpha$. 

By the assumption that $P^\alpha$ is $\varepsilon$-wqo, we can find $s\ltr t$ in $C$ such that 
$\varphi(h(s))\le_{P^{\alpha}} \varphi(h(t))$. Let $\rho:\,\alpha\to\alpha$ be a strictly increasing function witnessing that $\varphi(h(s))\le_{P^{\alpha}} \varphi(h(t))$.

We claim that if $\xi<\gamma_{h(s)}$ then $\rho(\xi)<\gamma_{h(t)}$. If this is not the case
then $\gamma_{h(t)}\leq \rho(\xi)$. Take $\zeta$ such that $\xi\leq \zeta<\gamma_{h(s)}$
and $h(s)(\zeta)\neq a$, which has to exist by the choice of $\gamma_{h(s)}$. Then $h(s)(\zeta)\leq_P 
(\varphi(h(t)))(\rho(\zeta))=a$. However, since $a$ is minimal, $h(s)(\zeta)= a$, which yields  a contradiction. Thus  the claim holds. 

Hence the restriction of $\rho$ to the
left interval of $h(s)$ maps it to the
left interval of $h(t)$. We would like to say that also the right interval of $h(s)$ is mapped by $\rho$ to the right interval of $h(t)$, but this not necessarily true.
But the order-type of the right interval of $h(s)$ is $\leq$ the
order-type of the right interval of $h(t)$, so we can change $\rho$ on the right interval of $h(s)$ to obtain
a strictly increasing $\rho'$ from $\dom(h(s))$ to $\dom(h(t))$, finally obtaining that $h(s)\leq_{P^{<\alpha}} h(t)$.

\end{proof}

\begin{lemma}\label{beta<alpha} (1) Suppose that $\beta<\alpha$ and that $P^{<\alpha}$ is
$\varepsilon$-wqo for some
$\varepsilon\in [\omega, \omega_1)$. Then so are $P^\beta$ and $P^{<\beta}$.

{\noindent (2)} Suppose that $P^{\alpha}$ is
$\varepsilon$-wqo for some
$\varepsilon\in [\omega, \omega_1)$. Then so is $P^{<\alpha}$.
\end{lemma}

\begin{proof} (1) The identity function is an order generating function from
$P^\beta$ to $P^{<\alpha}$ as well as from $P^{<\beta}$ to $P^{<\alpha}$. Hence the conclusion follows from Lemma \ref{orderpreserving}(2). 

{\noindent (2)} Similar, since the identity function is an order generating function from
$P^{<\alpha}$ to $P^{\alpha}$.
$\qed_{\ref{beta<alpha}}$
\end{proof}

\begin{lemma}\label{intoinitseg} Suppose that $\alpha>0$. Then there is an order-generating function of $P^\alpha$ into $I(P^{<\alpha})$.
\end{lemma}

\begin{proof}  To an $f\in P^\alpha$, we associate the initial segment generated by $\{f\rest\beta:\beta<\alpha\}$, that is 
\[
I(f)=\{s\in P^{<\alpha}:\,(\exists \beta<\alpha)\,s\le_{P^{<\alpha}} f\rest\beta\}.
\]
We claim that  $I(f)\subseteq I(g)\implies f\le_{P^\alpha} g$. To prove this, we shall define a strictly increasing function $\rho:\,\alpha\to\alpha$ such that for
all $\beta<\alpha$ we have $f(\beta)\le_P g(\rho(\beta))$. The definition of $\rho(\beta)$ is by induction on $\beta$. The inductive assumption is that 
for all $\beta<\alpha$, the value of $\rho(\beta)$ is the least possible such that 
\begin{itemize}
\item for all $\gamma<\beta$ we have $\rho(\gamma)<\rho(\beta)$ and
\item $f\rest \beta\le_{P^{<\alpha}} g\rest\rho(\beta)$.
\end{itemize}
At $\beta=0$ we have that $\emptyset=f\rest 0=g\rest 0$, so we can start with $\rho(0)=0$. Given $\beta>0$ such that $\rho(\gamma)<\alpha$ has been defined for all $\gamma<\beta$, let $\rho^\ast(\beta)=\sup\{\rho(\gamma):\,
\gamma<\beta\}$. 
If $\beta<\alpha$ is a limit, we have that 
\[
f\rest \beta\le_{P^{<\alpha}} g\rest\rho^\ast(\beta),
\]
by the inductive assumption, and clearly $\rho^\ast(\beta)$ is the least possible value which satisfies both properties required of $\rho(\beta)$. So, we would like to set $\rho(\beta)=\rho^\ast(\beta)$, but we must first assure ourselves that $\rho^\ast(\beta)<\alpha$. By the assumption $I(f)\subseteq I(g)$, so
we have that $f\rest\beta\in I(g)$ and in particular there is $o<\alpha$ such that $f\rest\beta\le_{P^{<\alpha}} g\rest o$. Let $o^\ast$ be the least
such $o$. Then by the assumption on the minimality of $f(\gamma)$ for $\gamma<\beta$, we have chosen $f(\gamma)\le o^\ast$ for any such $\gamma$, and 
hence $\rho(\beta)=\rho^\ast(\beta)\le o^\ast<\alpha$.
Hence, $\rho(\beta)$ is well defined and the induction continues.

If $\beta=\gamma+1$, then $\rho^\ast(\beta)=\rho(\gamma)<\alpha$. As above, there is the minimal $o<\alpha$ such that $f\rest\beta\le_{P^{<\alpha}} g\rest o$. By the inductive assumption $\rho(\gamma)<o$ for such $o$. It suffices to let $\rho(\beta)=o$.
$\qed_{\ref{intoinitseg}}$
\end{proof}

Hence, by putting together Lemma \ref{intoinitseg} and Theorem \ref{alpha<alpha}, we have the following conclusion:

\begin{theorem}\label{newknowledge} Suppose that $\alpha>0$. If
 $I(P^{<\alpha})$ is $\varepsilon$-wqo for some
$\varepsilon\in [\omega, \omega_1)$, then so are $P^\alpha$ and $P^{<\alpha}$.
\end{theorem}

\section{Products}\label{sec:products} Using Theorem \ref{bqoproducts},
one can easily make the following observation

\begin{observation}\label{finiteproduct}  A finite product of $\sigma$-bqo orders is $\sigma$-bqo.
\end{observation}

\begin{proof} By an inductive argument, it suffices to prove that if $P$ and $Q$ are $\sigma$-bqo then so is $P\times Q$.
Let $P=\bigcup_{n<\omega} P_n$ and $Q=\bigcup_{n<\omega} Q_n$, where each $P_n$ and $Q_n$ are bqo. Then
$P\times Q=\bigcup_{n,m<\omega} P_n\times Q_m$, which is $\sigma$-bqo by Theorem \ref{bqoproducts}.
$\qed_{\ref{finiteproduct}}$
\end{proof}

We prove that, on the condition of preserving the property of wqo, the countable infinite product also preserves the property of being $\sigma$-bqo.

\begin{theorem}\label{products} Let $1\le\alpha < \omega_1$. Suppose that $P$ is $\sigma$-bqo and that 
$P^\alpha$ is wqo. Then the orders $P^{<\alpha}$, $I(P^{<\alpha})$ and $P^\alpha$ are all $\sigma$-bqo.
\end{theorem}

\begin{proof} 
The proof is by induction on $\alpha$. For $\alpha=1$, the assumption is that $P$ is $\sigma$-bqo.
The order $P^{<\alpha}$ is the singleton containing the empty sequence, clearly a $\sigma$-bqo and similarly, 
$I(P^{<1})$ contains two elements, so is $\sigma$-bqo.

Arriving at $\alpha$, the assumption is that $P^\alpha$ is wqo. We use Lemma \ref{beta<alpha} to conclude that for each $\beta<\alpha$ the order $P^\beta$ is wqo. By Theorem \ref{alpha<alpha}, for each $\beta\le\alpha$, 
the order $P^{<\beta}$ is wqo.
Therefore, by the induction hypothesis, for each $\beta<\alpha$ the order $P^\beta$ is $\sigma$-bqo.  It follows that $P^{<\alpha}=\bigcup_{\beta<\alpha} P^\beta$ is $\sigma$-bqo. 

Next, we shall prove that $I(P^{<\alpha})$ is $\sigma$-bqo.
We have shown that $P^{<\alpha}$ is wqo, so each  subset of $P^{<\alpha}$ has a finite set of minimal elements. Calling such minimal elements {\em minimals}, let us define
\[
Q_n=\{A\in I(P^{<\alpha}):\,P^{<\alpha}\setminus A\mbox{ has exactly } n\mbox{ minimals}\}.
\]
Hence 
$ I(P^{<\alpha})=\bigcup_{n<\omega} Q_n$ and it suffices to show that each $Q_n$
is $\sigma$-bqo. Let $n$ be arbitrary. We have proved that $P^{<\alpha}$ is $\sigma$-bqo, so by Observation \ref{finiteproduct}, it follows that $(P^{<\alpha})^n$ is as well. Let 
\[
A_n=\{\bar{s}\in (P^{<\alpha})^n:\,\bar{s}\mbox{ is an antichain}\}.
\]
By Corollary \ref{fromalphatobeta}(3) it follows that $A_n$ is $\sigma$-bqo. The fact that  $I(P^{<\alpha})$ is $\sigma$-bqo will follow by the following Lemma \ref{mins} and an application of Lemma \ref{orderpreserving}(1).

\begin{lemma}\label{mins} There is an order-preserving surjection from $(A_n, \le_{(P^{<\alpha})^n})$ to the order
$(Q_n, \subseteq)$.
\end{lemma}

\begin{proof} (of Lemma \ref{mins}) We define the surjection $F$ as follows.
For $\langle x_i:\,i< n\rangle$ in $A_n$, we let 
\[
F(\langle x_i:\,i< n\rangle)=\{y\in P^{<\alpha}:\,(\nexists i<n) (x_i\le_{P^{<\alpha}} y)\}.
\]
First let us notice that $F(\langle x_i:\,i< n\rangle)$ is indeed an initial segment of $P^{<\alpha}$, so 
$F(\langle x_i:\,i< n\rangle)\in I(P^{<\alpha})$. We have that $P^{<\alpha}\setminus F(\langle x_i:\,i< n\rangle)$
consists of exactly those elements $y$ of $P^{<\alpha}$ such that for some $i<n$ we have $x_i\le_{P^{<\alpha}} y$,
so, since $\langle x_i:\,i< n\rangle$ is an antichain, the minimal elements of 
$P^{<\alpha}\setminus F(\langle x_i:\,i< n\rangle)$ are exactly $\{x_i:\,i<n\}$. Hence $F(\langle x_i:\,i< n\rangle)\in Q_n$.

To see that $F$ is a surjection, suppose that $A\in Q_n$ and let $\{y_i:\,i<n\}$ be the minimal elements of
$P^{<\alpha}\setminus A$. In particular, $\{y_i:\,i<n\}$ form an antichain of length $n$, so for any
arrangement of $\{y_i:\,i<n\}$ into a sequence $\bar{s}$, we have that $F(\bar{s})=A$.

Let us check that $F$ is order-preserving. Suppose that $\langle x_i:\,i< n\rangle\le_{(P^\alpha)^n}
\langle y_i:\,i< n\rangle$ for two antichains $\langle x_i:\,i< n\rangle$ and $\langle y_i:\,i< n\rangle$. Therefore $x_i\le_{P^\alpha} y_i$
for each $i<n$ and, so, if $y\in P^{<\alpha}$ is such that $(\nexists i<n) (x_i\le_{P^{<\alpha}} y)$, then clearly
$(\nexists i<n) (y_i\le_{P^{<\alpha}} y)$. So $F(\langle x_i:\,i< n\rangle)\subseteq F(\langle y_i:\,i< n\rangle)$.
$\qed_{\ref{mins}}$
\end{proof}

Finally, we shall prove that $P^\alpha$ is $\sigma$-bqo.
If $\alpha=\beta+1$ for some ordinal $\beta$, let us construct an order generating function of $P^\alpha$ into
the product $P^{\beta} \times P$ (ordered coordinatwise). Namely, we let $G(f)=(f\rest\beta, f(\beta))$.
Suppose that $G(f)\le_{P^{\beta} \times  P} G(g)$, in particular there is an increasing
function $\rho:\,\beta \to\beta$ such that $f(\gamma)\le _P g(\rho(\gamma))$ for all $\gamma<\beta$ and $f(\beta)\le_P
g(\beta)$. Then the function $\rho':\alpha\to\alpha$ defined by $\rho'=\rho\cup\{(\beta,\beta)\}$ shows that
$f\le_{P^\alpha} g$. Since the inductive hypothesis implies that $P^{\beta}$ is $\sigma$-bqo and we have assumed that $P$   is $\sigma$-bqo, it follows by Observation \ref{finiteproduct} that $P^{\beta} \times P$ is $\sigma$-bqo. By Lemma
\ref{orderpreserving}(2) it follows that $P^\alpha$ is $\sigma$-bqo.

For $\alpha$ a limit ordinal the above proof will not work, but here we can use that $I(P^{<\alpha})$
is $\sigma$-bqo and apply Lemma \ref{orderpreserving}(1), since Lemma \ref{intoinitseg} gives us that there is an order-generating
function of $P^\alpha$ into $I(P^{<\alpha})$.
This concludes the proof.
$\qed_{\ref{products}}$
\end{proof}

\section{Lexicographic Sums}\label{sec:lexsums}
The following is well known (see for example Exercise 9.14(iv) in \cite{SimpsonNW}).

\begin{theorem}\label{sumbqo} A lexicographic sum of bqo orders over a bqo order is bqo.
\end{theorem}

An easy consequence is the following theorem which shows that the class of $\sigma$-bqo is closed under lexicographic sums.

\begin{theorem}\label{lexsums} Suppose that $Q$ and $\langle P_q:\,q\in Q\rangle$ are $\sigma$-bqo orders.
Then so is the lexicographic sum 
\[
\Sigma_{q\in Q} P_q.
\]
\end{theorem}

\begin{proof} Denote $P=\Sigma_{q\in Q} P_q$, Let $Q=\bigcup_{n<\omega} Q_n$, where each $Q_n$ is a bqo. 
Hence $P=\bigcup_{n < \omega} \Sigma_{q\in Q_n} P_q$, and therefore it suffices to show that each 
$\Sigma_{q\in Q_n} P_q$ is $\sigma$-bqo. Fix $n<\omega$.

For each $q\in Q_n$ let $P_q=\bigcup_{ m< \omega} B^q_m$ where each $B^q_m$ is a bqo. Then we have
\[
\Sigma_{q\in Q_n} P_q =\bigcup_{m<\omega}  \Sigma_{q\in Q_n}B^q_m,
\]
and by Theorem \ref{sumbqo}, each $\Sigma_{q\in Q_n}B^q_m$ is a bqo.
$\qed_{\ref{lexsums}}$
\end{proof}

\section{Dress-Schiffels Product}\label{sec:DS} 
We start by a review of the relevant definitions and known facts.

\begin{definition}\label{def:Dress-Schieffels}
Let $(I,\leq_I)$ be a poset and suppose that for each $i\in I$ we have a poset $(P_i, \leq_i)$  with the smallest
element $0_i$. 

By $\prod_{i\in I}P_i$ we denote the set of functions $f$ defined on $I$ and such that $f(i)\in P_i$ for
every $i\in I$. For $f\in \prod_{i\in I}P_i$, the {\em support of $f$} is the set $\supp(f)=\{ i\in I:\,f(i)\neq 0_i\}$.
For $f,g \in \prod_{i\in I}P_i$ we let $\D(f,g) \deq \{ i\in I: f(i)\neq g(i)\}$

The {\em Dress--Schiffels poset} $\otimes^{\DS}_{i\in I} P_i$ has as its universe the set
\[
\prod_{i\in I}P_i^{\mbox{fin}}\deq\{ f\in \prod_{i\in I}P_i: \supp(f) \mbox{ is finite}\},
\]
and is ordered by $f\leq_{\DS} g$ iff 
\begin{quotation}
\fbox{$f(i)<_i g(i)$ for every $i$ $I$-maximal in $\D(f,g)$.}
\end{quotation}
\end{definition}

Note that for $i\in \D(f,g)$, $f(i)\leq_i g(i)$ implies $f(i) <_i g(i)$. Also, $f<_{\DS} g$ iff $\D(f,g)\neq \emptyset$
and for every $i$ that is maximal in $\D(f,g)$, $f(i)<_i g(i)$.

\begin{lemma}\label{lem:DSpo} The relation
$\leq_{\DS}$ is a partial ordering on $\otimes^{\DS}_{i\in I} P_i$.
\end{lemma}

\begin{proof} We need to prove that $\leq_{\DS}$ reflexive, antisymmetric and transitive.

Suppose that $f\in \prod_{i\in I}P_i$. Then $f\leq_{\DS} f$ because $\D(f,f)= \emptyset$. Surely, if $\D(f,g)=\emptyset$ then $f=g$, and hence
$f\leq_{DS} g\leq_{DS} f$ implies that $f=g$. Assume now that $f\leq_{\DS} g \leq_{\DS} h$. In order to prove that
$f\leq_{\DS} h$, consider any maximal $i\in \D(f,h)$ and suppose for a contradiction that $f(i)\not <_i h(i)$. 
Then \[f(i)\not\leq_i g(i)\ \text{or } g(i)\not \leq_i h(i)\] or else $f(i)\leq_i h(i)$, which implies $f(i)<_ih(i)$ since $i\in \D(f,h)$.

In the case $f(i)\not\leq_i g(i)$ we have that $i\in\D(f,g)$. Since $f\leq_{\DS} g$, we have that $i$ is not maximal in $\D(f,g)$ and hence
there is $j>_I i$ maximal in $\D(f,g)$ and
such that $f(j)<_{j} g(j)$. Note that for any index above $j$, the functions $f,g,h$ have the same value (since $i$ is maximal
in $\Delta(f,h)$ and $j>_I i$ is maximal in $\Delta(f,g)$).
By the maximality of $i$ in $\D(f,h)$, we have $f(j)=h(j)$.
Thus $h(j)<_j g(j)$ and hence $j\in \D(g,h)$ but is not maximal. So, there is some $k>_I j$ maximal in $\D(g,h)$ and such that $g(k)<_k h(k)$. 
But this contradicts our remark that the three functions get the same value above $j$.  The same line of argument
is used in the second case, $g(i)\not \leq_i h(i)$. 
$\qed_{\ref{lem:DSpo}}$
\end{proof}

Two simple examples of Dress-Schiffels product ordering can be obtained in the following way
\begin{enumerate}
\item In case that $I$ is a linear ordering, $\otimes^{\DS}_{i\in I} P_i$ can be thought of as an antilexicographical
ordering in which comparing two words begins with the highest index and checks down until a first difference is found.

\item In case we only take functions in $\otimes^{\DS}_{i\in I} P_i$ whose support is a singleton, we get
the sum $\Sigma_{i\in I} P_i$.
\end{enumerate}

Dress and Schiffels \cite{DressSchiffels} proved that if $I$ and $P_i$ for every $i\in I$ are wqo posets, then so is $\otimes^{\DS}_{i\in I} P_i$. 
We prove a more general theorem using the combinatorics of barriers, namely :

\begin{theorem}\label{generalisedDressShiffels} If $I$ and $P_i$ for every $i\in I$ are bqo ($\alpha$-wqo for some $\alpha\in [\omega,\omega_1)$, $\sigma$-bqo) posets, then so is $\otimes^{\DS}_{i\in I} P_i$. 
\end{theorem}

As it is somewhat tehnical, the full proof of Theorem \ref{generalisedDressShiffels} is given in Section \ref{sec:proofs}. As an indication, here we give a combinatorial proof of the original theorem. This proof was given by
Pouzet in \cite{Mauricebook} and it serves as a model for the proof of the bqo part of Theorem \ref{generalisedDressShiffels}.

\begin{theorem}[Dress and Schiffels]
\label{TDS}
If $I$ and $P_i$ for every $i\in I$ are wqo posets, then so is $\otimes^{\DS}_{i\in I} P_i$. 
\end{theorem}

\begin{proof}
Suppose for a contradiction that a sequence $\langle f_i\mid i\in\omega \rangle$ in $\otimes^{\DS}_{i\in I} P_i$ is 
such that for every $m<n$ we have $f_m \not\leq_{\DS} f_n$. So there is an index $i=i_{m,n}\in I$ maximal in $\D(f_m,f_n)$ and
 such that $f_m(i)\not<_i f_n(i)$. So \[f_m(i)\not\leq_i f_n(i)\] since $i\in \D(f_m,f_n)$. In particular $f_m(i)>0_i$ 
(or else $f_m(i)\leq f_n(n)$). Thus for every $n>m$, $i_{m,n}\in\supp(f_m)$, which is a finite subset of $I$.
 Thus, for any $m\in \omega$ and any infinite set $X\subseteq \omega$
there is an infinite $X'\subseteq X$ (of numbers above $m$) such that for some fixed $i\in\supp(f_m)$, for all $n\in X'$,
$i=i_{m,n}$. With this observation in mind we can construct (in an inductive process) an infinite subset $X\subseteq\omega$ and assign to 
each $m\in X$ some
$i_m\in I$ such that for every $n\in X$ such that $m<n$, $i_m=i_{m,n}$ is maximal in $\D(f_m,f_n)$ and 
\begin{equation}
\label{Eq1}
f_m(i_m)\not \leq_{i_m} f_n(i_m).
\end{equation}
 Using the assumption that $I$ is wqo, there
is an infinite subset $X_1$ of $X$ such that for every $m<n$ both in $X_1$, $i_m\leq_I i_n$. Using Ramsey's theorem with view
of answering if $i_m=i_n$ or not, we may conclude that there is
an infinite subset $X_2\subseteq X_1$ which exhibits one of the two following behaviours.
In the first case, we have that for some fixed $i\in I$,  $i=i_n $ for every $n\in X_2$.  Then, by (\ref{Eq1}), the sequence
$f_n(i)$, for $n\in X_2$ is a bad sequence in $P_i$, which is impossible. 
In the second case, $i_m<_I i_n$ whenever $m<n$ are in $X_2$. Take a triple $m<n<k$ in $X_2$ such that $i_n\not\in \supp(f_m)$. Since  
$i_m\in \D(f_m,f_n)$ is maximal there, and as $i_m<_I i_n$, $f_m(i_n)=f_n(i_n)$. 
 But $f_m(i_n)=0_{i_n}$, and
hence $f_n(i_n)=0_{i_n}\leq f_k(i_n)$, which contradicts the definition of $i_n= i_{n,k}$. 
$\qed_{\ref{TDS}}$
\end{proof}

\section{$\alpha$-FAC Orders from $\alpha$-wqo Spines}\label{sec:Hausdorff} As mentioned in Section \ref{sec:present_results}, a poset $P$ is said to be {\em FAC} if all its antichains are finite. 

One of the main results of \cite{AbBo} is that all FAC orders can be built from wqo orders in a Hausdorff-style hierarchy. We shall show
an analogous result for $\alpha$-wqo orders and $\alpha$-FAC orders 
The following Lemma summarises Lemma 3.1 in \cite{AbBo}
(see also \cite[\S 9.9.1]{Fraissebook}) together with the discussion immediately following that
lemma. The lemma and a proof sketch can be found as Lemma 2.13 in \cite{abcdzt}. The following definition is used.

\begin{definition}\label{def:augmentation} A poset $(Q, \le_Q)$ is {\em an augmentation} of a poset $(P,\le P)$, if $Q$ has the same domain as $P$
and $\le_P\subseteq \le Q$ (that is, new relations may be added in $Q$). 
\end{definition}
    
\begin{lemma}[Abraham and Bonnet]\label{supertechnical} 
Let $P$ be a FAC poset. Then there exist an ordinal
$\zeta$ and $p_\beta \in P$, $Z_\beta \subseteq P$ for $\beta < \zeta$ such that:

\begin{enumerate}

\item  $\{ p_\beta : \beta < \zeta \}$ is a wqo subset of $P$. 

\item  $Z_\beta$ is a convex set with maximum element $p_\beta$. The $Z_\beta$ are disjoint and
      form a partition of $P$. 

\item   $P$ is an augmentation of the lexicographic sum of the
    $Z_\beta$ along the index set $(\zeta, \preceq)$, where 
    $\beta \preceq \gamma \iff p_\beta \le p_\gamma$. 
\end{enumerate}

\end{lemma} 

In the context of Lemma \ref{supertechnical}, call the order $\{ p_\beta : \beta < \zeta \}$ {\em a spine} of $P$.
Then we define:

\begin{definition}\label{allphafac} A FAC order is {\em $\alpha$-FAC} if it has an $\alpha$-wqo spine.
\end{definition}

We obtain a classification theorem about  $\alpha$-FAC orders following the lines of Section 6
in \cite{abcdzt}. The idea is, given
 a class $\mathcal{G}_0$ of linear orderings, we consider the class
  $\mathcal{G}$ of $\alpha$-FAC posets in which every chain belongs to $\mathcal{G}_0$. We prove
  that under mild assumptions on $\mathcal{G}_0$, every element of $\mathcal{G}$ can be
 built up from members of $\mathcal{G}_0$.

\begin{definition}\label{reasonable} A class $\mathcal{G}_0$ of linear orderings is {\em reasonable}
    if and only if $\mathcal{G}_0$ contains a nonempty ordering and
   is closed under reversals and restrictions.
\end{definition}

\begin{definition}\label{reasonable2} Given a reasonable class $\mathcal{G}_0$ of linear orderings,
    the $\alpha$-{\em closure} $\alpha cl( {\mathcal G}_0 )$ of $\mathcal{G}_0$ is the least
   class of posets which contains $\mathcal{G}_0$ and is closed
   under the operations:
\begin{itemize}
\item Lexicographic sum with index set  either an  $\alpha$-wqo poset, the inverse
   of a  $\alpha$-wqo poset, or an element of $\mathcal{G}_0$.  
\item Augmentation.
\end{itemize}
\end{definition} 

 It is easy to see that $\alpha cl( {\mathcal G}_0 )$ consists of FAC posets and is
 closed under restrictions and reversals. 
 
 The following technical lemma records another useful closure property of 
$\alpha cl( {\mathcal G}_0 )$. Here $q^\parallel$ is $\{r\in Q:\,q, r\ \mbox{comparable}\}$ and 
$q^\perp$ is $\{r\in Q:\,q, r\ \mbox{incomparable}\}$.

\begin{lemma}  \label{technical} 
 Let $Q$ be a poset, and suppose that there is  $q \in Q$ be such that
 $q^\parallel$ and $q^\perp$ are both in
 $\alpha cl( {\mathcal G}_0 )$. Then $Q \in \alpha cl( {\mathcal G}_0 )$.
\end{lemma} 

\begin{proof}  Let $X_0 = q^\parallel$ and
   $X_1 = q^\perp$. 
   Now we form a lexicographic sum $Q'$
   of the orderings $X_i$, ordering the indices so that 
    $0$ is incomparable with $1$. Since the $X_i$ both lie
   in $\alpha cl( {\mathcal G}_0 )$ and finite posets are trivially $\alpha$-wqo,
   $Q' \in \alpha cl( {\mathcal G}_0 )$. Clearly $Q$ is an augmentation
   of $Q'$, so $Q \in \alpha cl( {\mathcal G}_0 )$. 
\end{proof}

 The main theorem of this section is Theorem \ref{forbiddencor}, which tells us how to obtain exactly the class of $\alpha$-wqos 
 as the $\alpha$-closure of a certain reasonable class of $\alpha$-FAC orders. On the way to that theorem, we prove Theorem \ref{forbiddenlo} 
which is a general statement from which we obtain Theorem \ref{forbiddencor} as a consequence.
 
\begin{theorem} \label{forbiddenlo} Let $\mathcal{G}_0$ be a reasonable class of posets and
 let $\mathcal{G}$ be the class of $\alpha$-FAC posets such that every chain
  is in $\mathcal{G}_0$. Then  $\mathcal{G} \subseteq \alpha cl( {\mathcal G}_0 )$.
\end{theorem}

\begin{proof}   Note that any FAC poset $P$ can be equipped with a rank function, defined by considering the
order of reverse inclusion on the set of all antichains. Since the antichains are finite, this order is well-founded and hence allows for
a rank function, that is, an order-preserving function into an ordinal. The smallest such ordinal is called the antichain rank of $P$.

We will proceed by induction on the antichain rank $\rho(P)$ of an $\alpha$-FAC poset
    $P \in \mathcal{G}$.  We note that since $\mathcal{G}_0$ is closed under
    restriction, $\mathcal{G}_0 \subseteq \mathcal{G}$.

     Let $P$ be an $\alpha$-FAC poset, and assume that 
    $Q \in \alpha cl( {\mathcal G}_0 )$ for every $\alpha$-FAC order $Q \in \mathcal{G}$ with
    $\rho(Q) < \rho(P)$. In particular $p^\perp \in \alpha cl( {\mathcal G}_0 )$ for 
    every $p \in P$, a fact which will play a crucial role at several 
    points.

    We define a binary relation $\equiv$ on $P$ by stipulating that 
    $p \equiv q$ if and only if:
\begin{itemize}
\item $p$ is incomparable with $q$, or
\item $p \le q$ and $(p, q) \in \alpha cl( {\mathcal G}_0 )$, or
\item $q \le p$ and $(q, p) \in \alpha cl( {\mathcal G}_0 )$.
\end{itemize}

\begin{claim} $\equiv$ is an equivalence relation.
\end{claim} 

 Clearly $\equiv$  is reflexive and symmetric, so we check only that it is
     transitive. Let $a \equiv b$ and $b \equiv c$,
     where we may as well assume that $a, b, c$ are distinct. 
     There are four cases to check:

\medskip

\noindent Case 1:  $a < b$ and $b < c$. Let $Q = (a, c)$. Then
   $(a, b)$ and $(b, c)$ are in $\alpha cl( {\mathcal G}_0 )$
   by definition. $b^\perp \in \alpha cl( {\mathcal G}_0 )$ by our assumption
   on $P$, and $b^\perp \cap Q \in \alpha cl( {\mathcal G}_0 )$ 
   because $\alpha cl( {\mathcal G}_0 )$ is closed under restriction. 
   Applying Lemma \ref{technical} we see that
   $Q \in\alpha cl( {\mathcal G}_0 )$, and so by definition
   $a \equiv c$. 

\medskip

\noindent Case 2: $a < b$ and $b \perp c$. If $a \perp c$ we are done, so we
   assume that $a < c$. Let $X_0 = (a, c) \cap (a,b)$, and 
   $X_1 =  (a, c) \setminus (a,b)$. Then $X_0 \in \alpha cl( {\mathcal G}_0 )$
   because $(a, b) \in \alpha cl( {\mathcal G}_0 )$, and  
   $X_1 \in \alpha cl( {\mathcal G}_0 )$ because $X_1 \subseteq b^\perp$.
   Finally $(a, c)$ is an augmentation of the lexicographic sum
   of $X_0, X_1$ along $\{0,1\}$ in which $0, 1$ are incomparable.

\medskip

\noindent Case 3: $a \perp b$ and $b < c$. This is exactly like the previous case.

\medskip

\noindent Case 4: $a \perp b$ and $b \perp c$. If $a \perp c$ we are done, so
   we may as well assume that $a < c$. Then $(a, c) \subseteq b^\perp$,
   and again we are done.

\medskip

\begin{claim} The equivalence classes of $\equiv$ are convex.
\end{claim}

   Let $a < b$ with $a \equiv b$. Then for every $c \in (a, b)$
   we have $(a, c) \subseteq (a, b)$, and so $a \equiv c$ since
   $\alpha cl( {\mathcal G}_0 )$ is closed under restriction.

\begin{claim} Each equivalence class is in $\alpha cl( {\mathcal G}_0 )$.    
\end{claim} 

  Let $C$ be such a class, and let $c \in C$. By Lemma \ref{technical}
  and  the fact that $c^\perp \in \alpha cl( {\mathcal G}_0 )$,  it will suffice to show
  that $\{ d \in C: d > c \}$ and $\{ d \in C: d < c \}$ 
  are both in  $\alpha cl( {\mathcal G}_0 )$.
    
   In fact we will just argue that the set  $Y = \{ d \in C: d > c \}$ 
   is in   $\alpha cl( {\mathcal G}_0 )$, the argument for 
   $\{ d \in C: d < c \}$ will be symmetric.
   Since $Y$ is an $\alpha$-FAC poset we may
   appeal to Lemma \ref{supertechnical} and choose an ordinal $\zeta$ and
   $d_\beta \in Y$, $Z_\beta \subseteq (c, d_\beta)_P \subseteq Y$ for $\beta < \zeta$
   such that
\begin{enumerate}
\item $\{ d_\beta :\beta < \zeta \}$ is a $\alpha$-wqo subset of $Y$.
\item Each $Z_\beta$ is convex with maximum element $d_\beta$, and the $Z_\beta$ form a 
    partition of $Y$.  
\item  $Y$ is an augmentation of the lexicographic sum
   of the $Z_\beta$ for $\beta < \zeta$, with indices ordered
   by $\beta \prec \gamma \iff d_\beta < d_\gamma$.
\end{enumerate}

 Since $C$ is an equivalence class,
   $(c, d_\beta)_P \in \alpha cl( {\mathcal G}_0 )$, and so $Z_\beta \in \alpha cl( {\mathcal G}_0 )$.
   It follows from the closure properties of $\alpha cl( {\mathcal G}_0 )$
   that $Y \in \alpha cl( {\mathcal G}_0 )$.

\begin{claim} If $C$ and $D$ are distinct equivalence classes,
   then either every element of $C$ is less than every element of $D$
   or vice versa.
\end{claim} 

Since incomparable elements are equivalent,
   every $c \in C$ is comparable with every $d \in D$.
 Suppose for a contradiction that  we have
   $c_1 < d < c_2$ with $c_1, c_2 \in C$ and $d \in D$.
 Since classes are convex we have $d \in C$, a contradiction since
 $C$ and $D$ are disjoint.

\begin{claim} $P \in  \alpha cl(\mathcal{G}_0)$.
\end{claim}

   The equivalence classes are linearly ordered in some order type $L$,
   and $P$ is an $L$-indexed sum of equivalence classes. In particular
   (choosing a point from each class) $P$ contains a copy of $L$.
   Since every chain of $P$ lies in $\mathcal{G}_0$, $L \in \mathcal{G}_0$.
   Finally since each class is in 
   $\alpha cl(\mathcal{G}_0)$, and  $\alpha cl(\mathcal{G}_0)$ is closed under lexicographic
   sums with index set in $\mathcal{G}_0$, we see that
   $P \in  \alpha cl(\mathcal{G}_0)$.
    $\qed_{\ref{forbiddenlo}}$
\end{proof}

    Theorem \ref{forbiddenlo} tells us that every element of $\mathcal{G}$
    is built up from elements of $\mathcal{G}_0$ by  a certain recipe.
    It does {\em not} in general guarantee that all posets built up
    according to this recipe are in $\mathcal{G}$, because 
    for example
    $\mathcal{G}$ may not be closed under augmentations. However, under an additional assumption we obtain

\begin{theorem}
\label{forbiddencor} Let $\mathcal{G}_0$ be a reasonable class of posets and
 let $\mathcal{G}$ be the class of $\alpha$-FAC posets such that every chain
  is in $\mathcal{G}_0$
    and assume in addition that:
\begin{enumerate}
\item $\mathcal{G}_0$ contains all well-orderings, and is closed under
   lexicographic sums with index set in $\mathcal{G}_0$.
\item $\mathcal{G}$ is closed under augmentations.
\end{enumerate}

   Then $\mathcal{G} = \alpha cl(\mathcal{G}_0)$.

\end{theorem}

\begin{proof} The extra closure assumptions on $\mathcal{G}_0$
   easily imply that $\mathcal{G}$ is closed under lexicographic sums
   with index sets that are $\alpha$-wqo, converse $\alpha$-wqos or elements of
   $\mathcal G$. Since we  also assumed that $\mathcal{G}$ is
   closed under augmentation, $\mathcal{G}$ is closed under 
   all the operations which are used to build 
   $\alpha cl(\mathcal{G}_0)$. Since $\mathcal{G}_0 \subseteq \mathcal{G}$
   it follows that $\alpha cl(\mathcal{G}_0) \subseteq \mathcal{G}$,
   and hence by Theorem \ref{forbiddenlo} that  
   $\alpha cl(\mathcal{G}_0) = \mathcal{G}$.
   $\qed_{\ref{forbiddencor}}$
\end{proof}

\begin{example}\label{} Taking $\mathcal{G}_0$ to be the class of all linear orders, we obtain that the class of
all $\alpha$-FAC orders is the $\alpha$-closure of the class of all linear orders. That is, $\alpha$-FAC orders
are exactly the closure of the class of linear orders under lexicographic sums along $\alpha$-wqos orders and augmentations.

In particular, with $\alpha=\omega^2$, we obtain that the closure of the class of linear orders under lexicographic sums along wqo's which do not contain a copy of Rado's
order and augmentations gives all $\omega^2$-FAC orders.
\end{example}

\section{Some More Technical Proofs}\label{sec:proofs} In this final section we give various proofs that were mentioned in the main body of the paper as being here. We also give a few more technical notions about barriers that were used in our proofs. 

Given a barrier $U$, a {\em successive sequence} of members of $U$ is a sequence of elements of $U$ such that each member of the sequence
stands in the $\triangleleft$ relation with the next one. The following lemma is known and easily proved.

\begin{lemma}\label{lem:succesive} Let $U$ be a barrier and let $r,s\in U$ be such that all members of $r$ are strictly
smaller than all members of $s$ (which can be written as $r\ll s$). Then there exists a unique successive sequence of members of $U$ of length $|r|+1$ that
starts with $r$ and ends with $s$. If $r= \{a_1,\ldots,a_m\}$ is an increasing enumeration of $r$, then this successive
sequence has length $m+1$, $r_1,\ldots,r_{m+1}$, and is such that 
 \[r=r_1\ltr r_2\ltr\cdots\ltr r_{m+1}=s, \] 
and  $a_i=\min (r_i)$ for $1\leq i \leq m$, and $r_i\subseteq r\cup s$ for $1\leq i \leq m+1$.
 (It follows that $|r|\leq |s|$ whenever $r\ll s$ are in $U$.)
\end{lemma}

For example, if $r=\{3,7\}$ and $s=\{16,20\}$, then the sequence $\{3,7\}$, $\{7,16\}$, $\{16,20\}$ is a successive sequence.

\begin{corollary}\label{cor:leading}
Let $r,s$ be any two members of a barrier $U$. Then there exist two successive sequences
in $U$, one starting with $r$ and the other with $s$, that have the same last element. 
\end{corollary}

For a proof take a member $t$
of $U$ that lies above both $r$ and $s$  (that such a $t$ exists from the property (2)(2) in the definition of a barrier, Definition \ref{def: barrier} )
and apply Lemma \ref{lem:succesive} to find a successive sequence that leads from $r$ to $t$ and a successive sequence from $s$ to $t$.

\begin{lemma}
\label{Lem3}
Suppose that $D$ is a barrier and $f$ is a function defined over $D$ such that whenever $d_1,d_2\in D$ and $d_1\ltr d_2$ then $f(d_1)=f(d_2)$.
Then $f$ is constant on $D$: for every $d_1,d_2\in D$, $f(d_1)=f(d_2)$.
\end{lemma}

\begin{proof} Let $d_1,d_2\in D$, and find a successive sequence in $D$ that leads from $d_1,d_2\in D$, as per Corollary \ref{cor:leading}. Then necessarily $f(d_1)=f(d_2)$.
$\qed_{\ref{Lem3}}$
\end{proof}

To get a feeling for the following Lemma \ref{Lem5} consider a particularly simple case which was already used in the proof of Theorem \ref{TDS}.
 Let $B$ be the barrier of all singletons
$B=\{\, \{m\}\mid m\in \omega\}$, and then $B^2$ is the barrier of all pairs $\langle m,n\rangle$ such that $m<n$ in $\omega$.
Let $S$ be an arbitrary set and $g:B^2\to S$ an arbitrary function into $S$. Suppose that $\beta:B\to [S]^{<\omega}$
is an ``approximation'' function, in the sense that if $a,b\in B$ and $a\ltr b$, then $g(a\cup b)\in \beta(a)$. That is,
$\beta(a)$ gives a finite range of possible values of $g(\{a,b\})$ for $a<b$. Then for any $m\in\omega$ and infinite
set $X\subseteq \omega$ of numbers above $m$, there is an infinite subset $X'\subseteq X$ such that for every $n_1,n_2\in X'$,
$g(\{m,n_1\}) = g(\{m,n_2\})$. Using this observation in an inductive construction of a descending sequence of infinite
subsets of $\omega$ we obtain an infinite set $Y\subseteq\omega$ and for every $m\in Y$ a value $\ii(m)\in S$ exists such that
for every $n\in Y$ above $m$, $g(\{m,n\})=\ii (m)$. The following `choice' Lemma \ref{Lem5}  generalises these observations.

\begin{lemma}
\label{Lem5}
Suppose that $B$ is a barrier, $I$ an arbitrary set, and  $g:B^2\to I$ a  function 
such that for some function $\beta:B\to [I]^{<\omega}$
\[ \forall [a,b]\in B^2: g([a,b])\in \beta(a).\]
 Then there are a sub-barrier
$C\subseteq B$ and a function $\ii:C\to I$, such that for every $[a,b]\in C^2$,
$g([a,b])=\ii(a)$.
\end{lemma}

\begin{proof}
Say $B$ is a barrier and suppose for notational simplicity that $\omega=\bigcup B$ is its base.
 Define by induction on $n\in\omega$ a pair $X_n,Y_n$ and a function $\ii:B\cap \Pset(X_n)\to I$ such that:

\begin{enumerate}
\item
$\langle X_n:\,n<\omega\rangle$ is an increasing sequence of finite sets and $\langle Y_n:\,n<\omega\rangle$ is a decreasing sequence of
infinite sets.
\item $X_0=\emptyset$ and $Y_0=\omega$. 
\item $X_n \ll Y_n$. 
For all $a\in B\cap \Pset(X_n)$,  $\ii(a)\in \beta(a)$ and, for all $b\in B \cap [X_n\cup Y_n]^{<\omega}$,
$a\ltr b\Rightarrow g([a,b])= \ii(a).$
\end{enumerate}

It should be clear that this construction, if successful, gives the required barrier over $X=\bigcup_{n\in\omega} X_n$, namely $C= B\cap [X]^{<\omega}$, and $\ii(a)$ is defined for all $a\in B\cap [X]^{<\omega}$.

For the inductive definition, suppose that $(X_n,Y_n)$ is defined and has the required properties.
 Define $X_{n+1} = X_n\cup \{ p\}$ where  $p= \min(Y_n)$. Define 
$Y_{n+1}\subseteq Y_n\setminus \{ p\}$
by the following procedure. Enumerate the set $E=\{ a\subseteq X_{n+1} \mid a\in B \wedge \max(a)= p \}$  as $E= \{a_i\mid i<m\}$, and define a descending sequence of infinite subsets
$Z^0\supseteq Z^1\supseteq \cdots\subseteq Z^m$, beginning with $Z^0 = Y_n\setminus \{ p \}$, and with the aim of setting $Y_{n+1}=Z^m$. Having defined $Z^i$ we  define
$Z^{i+1}\subseteq Z^i$ in order to ensure item (3) with respect to $a=a_i$. In the following we use the notation $_*a$ for $a\setminus \{\min(a)\}$.

Consider $B' = \{b\subset Z^i\mid\ _*a\cup b\in B\}$. It is a barrier with base $Z^i$, since if $S\subseteq Z^i$
is an infinite subset, then $_* a \cup S$ is a subset of $X_n\cup Y_n$ which has an initial segment $x\in B$. Since $x$
is not an initial segment of $_*a$ (for $B$ is an inclusion antichain), the set $b=x\setminus\, _*a$ is in $B'$.
If $b\in B'$ then $a\cup b\in B^2$, and we partition $B'$ in accordance with the value of $g(a\cup b)\in \beta(a)$. 
Since $\beta(a)$ is finite, Nash-Williams' theorem can be used, and there is a fixed value $k_a\in \beta(a)$ 
and a sub-barrier $B''$ of $B'$ such that
for any $b\in B''$,
$g(a\cup b)=k_a$.  We define $Z^{i+1}$ to be the base of $B^{''}$, and $\ii(a)=k_a$. Property (3) holds for
$(X_{n+1},Y_{n+1})$. To see this, take any $a\subseteq X_{n+1}$ in $B$; there are two cases: $p\in a$ and $p\not\in a$. 
Suppose first that $p\in a$. 
Suppose that $b\in B$, is such that $a\ltr b$ and 
$b\subset X_{n+1}\cup Y_{n+1}$. For some $i<m$ the definition of $Z^{i+1}$
was made to ensure that $g(a\cup b)=k_a$. Now suppose that $p\not\in a$.
Then $a\subseteq X_n$, and if $b\in B$ is such that 
$b\subseteq X_{n+1}\cup Y_{n+1}\subseteq X_n\cup Y_n$ then
the induction hypothesis ensures that $g(a\cup b)=\ii(a)$. 
$\qed_{\ref{Lem5}}$
\end{proof}

\begin{proof-} ({\bf of Theorem \ref{generalisedDressShiffels}.}) We divide the proof in three parts.

 \noindent{(a)} Let $\alpha\in[\omega,\omega_1)$ be such that  $I$ and every $P_i$ for $i\in I$ 
are $\alpha$-wqo posets.
Let $B$ be a barrier of order-type $\alpha$ and let $f:B\to \otimes^{\DS}_{i\in I} P_i$ be some supposedly bad barrier sequence. We shall derive a contradiction.

For $b\in B$, we prefer to write $f_b$ rather than $f(b)$ for the value of the sequence $f$ at $b$. Thus for $i\in I$, $f(b)(i)\in P_i$ will be written $f_b(i)$.
Since $f$ is bad, for every $a,b\in B$ such that $a\ltr b$, $f_a\not\leq_{\DS} f_b$, and this means that there is an index $i=i_{a,b}$ maximal
in $\D(f_a,f_b)$ and such that 
\begin{equation}
\label{E1a}
f_a(i)\not \leq_{P_i} f_b(i).
\end{equation}
Note that $i_{a,b}\in \supp(f_a)$ which is a finite set, and if we define
$\beta(a)= \supp(f_a)$, then we are in the framework of Lemma
 \ref{Lem5} (taking $i_{a,b}$ for $g([a,b]))$. This lemma gives a sub-barrier
$C\subseteq B$ and a function $\ii:C\to I$, such that
\begin{equation}
\label{Eq6}
\forall a,b\in C\, (a\ltr b \Rightarrow i_{a,b}=\ii(a)).
\end{equation}
 Thus, for $a,b\in C$,
\begin{equation}\label{E1}
a\ltr b \Rightarrow f_a({\ii}(a))\not\leq_{P_{{\ii}(a)}} f_b({\ii}(a)).
\end{equation}

 Since $I$ is $\alpha$-bqo, and $\ii:C\to I$ is
a barrier sequence, there exists a 
sub-barrier $D\subseteq C$ on which $\ii$ is perfect:
 for every $a,b\in D$,
\begin{equation}
a\ltr b \Rightarrow \ii(a) \leq_I \ii(b).
\end{equation}

We partition the barrier $D^2$ into two classes. In the first class we put all $[a,b]\in D^2$
for which $\ii(a)=\ii(b)$, and in the second those for which $\ii(a)<_I \ii(b)$. Applying the Nash-Williams
partition theorem to $D^2$ we get a sub-barrier $E\subseteq D$ such that all members of $E^2$ are in the same
class. Now we distinguish two possible cases.

Case 1: for all $[a,b],[b,c]\in E^2$, $i_{a,b}=i_{b,c}$. So, by Lemma \ref{Lem3},
for some fixed $i\in I$, for every $[a,b]\in E$, $i=i_{a,b}$.
This implies that if $a,b\in E$ are such that $a\ltr b$ then $f_a(i)\not\leq_{P_i} f_b(i)$ (by (\ref{E1})).
However, this is in contradiction to our assumption about $P_i$ being $\alpha$-wqo.

Case 2: for all $[a,b],[b,c]\in E^2$, $i_{a,b} <_I i_{b,c}$. Pick a successive sequence of sets in $E$,   
 $a_0\ltr a_1\ltr a_2\ltr\cdots \ltr a_{m-1}\ltr a_m$, for $m$ large enough  such that $|\supp(f_{a_0})|<m$. 
So $i_{a_0,a_1}<_I i_{a_1,a_2}<_I i_{a_2,a_3}<_I \cdots<i_{a_{m-1},a_m}$ are all different (and $i_{a_0,a_1}\in \supp(f_{a_0})$).
Hence for some index $n$, $0<n < m$, $i_{a_{n}, a_{n+1}}\not\in \supp(f_{a_0})$. 
So
\[ f_{a_0}( i_{a_{n},a_{n+1}})= 0_{i_{a_{n}}}.\]
Yet, for every $k$, $i_{a_k,a_{k+1} }$ is maximal in $\Delta(f_{a_k},f_{a_{k+1}})$, and hence
 $f_{a_k}(j) = f_{a_{k+1}}(j)$ whenever $i_{a_k,a_{k+1}}<_I j$. 
Therefore, for $j=i_{a_{n},a_{n+1} }$,
 \[f_{a_0}(j) = f_{a_1} (j)=\cdots= f_{a_{n-1}}(j)=f_{a_{n}}(j).\]
So $f_{a_{n}}(j)=0_{j}$. But this contradicts $a_{n}\ltr a_{n+1}$ and 
$f_{a_{n}}({j} )\not\leq_j f_{a_{n+1}}( j)$.

\medskip {\noindent (b)} The case of bqo. This case follows from (a) since by Definition \ref{nashwilliams}, a poset is bqo iff it is $\alpha$-wqo for every $\alpha<\omega_1$.

\medskip {\noindent (c)} In this case we deal with $\sigma$-bqo orders. Let $I=\bigcup_{n<\omega} I_n$ and for each $i\in I$, let
$P_i= \bigcup_{n<\omega} P^n_i$, where each $I_n$ and $P^n_i$ are bqo. Without loss of generality, each of these unions is increasing, since it is easily seen that the finite union of bqos is still a bqo. By the definition of what is meant by the union of orders, each union is made so that each $I_n$ is an induced suborder of $I$ (i.e. no new relations between elements of $I_n$ are added when passing to $I$) and similarly for each $P^n_i$. We shall now prove a series of lemmas.

\begin{lemma}\label{lem:prvared} $\prod_{i\in I}P_i^{\rm fin}= \bigcup_{n<\omega} \prod_{i\in I}(P^n_i)^{\rm fin}$.
\end{lemma}

\begin{proof} Let $a\in \prod_{i\in I}P_i^{\rm fin}$, hence $\supp(a)$ is some finite set $F$. For each $i\in F$ let $n_i$ be such that $a(i)\in P_i^{n_i}$
and let $n^\ast=\max\{n_i:\,i\in F\}$. Then $a\in \prod_{i\in I}(P^i_n)^{\rm fin}$. We have shown that $\prod_{i\in I}P_i^{\rm fin}\subseteq
\bigcup_{n<\omega} \prod_{i\in I}(P^n_i)^{\rm fin}$, and the other direction is trivially true. Hence the sets in the announcement of the lemma are the same.
$\qed_{\ref{lem:prvared}}$
\end{proof}

\begin{lemma}\label{lem:drugared} If $a,b\in \prod_{i\in I}P_i^{\rm fin}$ then $a\le_{DS}b$ with respect to the product
$\prod_{i\in I}P_i^{\rm fin}$  iff $a\le_{DS}b$ with respect to the product
$\prod_{i\in I}(P_i^n)^{\rm fin}$ for some $n$.
\end{lemma}

\begin{proof} Let $n$ be such that both $a,b\in \prod_{i\in I}(P_i^{n})^{\rm fin}$, which exists by Lemma \ref{lem:prvared}.
Then we have $a\le_{DS}b$ with respect to the product $\prod_{i\in I}P_i^{\rm fin}$  iff $a\le_{DS}b$ with respect to the product
$\prod_{i\in I}(P_i^{n^\ast})^{\rm fin}$.
$\qed_{\ref{lem:drugared}}$
\end{proof}

We conclude that $\otimes^{\DS}_{i\in I} P_i=\bigcup_{n<\omega} \otimes^{\DS}_{i\in I}(P^i_n)$. Hence it suffices to show that each 
$\otimes^{\DS}_{i\in I}(P^i_n)^{\rm fin}$ is a $\sigma$-bqo.

\begin{lemma}\label{lem:trecared} Let $I$ be as above and suppose that each $Q_i$ for $i\in I$ is a bqo, for some orders $Q_i$. Then $\otimes^{\DS}_{i\in I}(Q^i)$ is $\sigma$-bqo.
\end{lemma}

\begin{proof} 
By part (b) of the proof of Theorem \ref{generalisedDressShiffels}, which we completed, we have that for every $n\in\omega$, the product
$\otimes^{\DS}_{i\in I_n}(Q^i)$ is bqo. We shall produce an order-generating function from $\otimes^{\DS}_{i\in I}(Q^i)$ into 
$\bigsqcup_{n<\omega}\otimes^{\DS}_{i\in I_n}(Q^i)$. This notation means that we are taking a disjoint union, in which for $n\neq m$
there are no relations between the elements in $\otimes^{\DS}_{i\in I_n}(Q^i)$ and $\otimes^{\DS}_{i\in I_m}(Q^i)$ and in the copy of
$\otimes^{\DS}_{i\in I_m}(Q^i)$ the order is inherited from $\otimes^{\DS}_{i\in I_m}(Q^i)$. The proof will then follow by Lemma \ref{orderpreserving}(2), since $\bigsqcup_{n<\omega}\otimes^{\DS}_{i\in I_n}(Q^i)$ is $\sigma$-bqo. Let us denote by $\le^*$ the order in $\bigsqcup_{n<\omega}\otimes^{\DS}_{i\in I_n}(Q^i)$.

We shall use the notation $a\rest I_n$ for the restriction of $a$ to the coordinates in $I_n$, for $a\in \otimes^{\DS}_{i\in I}(Q^i)$. For each such $a$,
let $n^\ast_a$ be the minimal $n$ such that $\supp(a)\subseteq n$ and let us define $g(a)=a\rest n^\ast_a$. Suppose then that in 
$\bigsqcup_{n<\omega}\otimes^{\DS}_{i\in I_n}(Q^i)$ we have $g(a)\le^\ast g(b)$. This in particular means that $n^\ast_a=n^\ast_b$ (denote this value by $n$)
and so $\Delta(a,b)=\Delta(a\rest I_n, b\rest I_n)$. Now we notice that the $I$-maximal elements in $\Delta(a,b)$ are the same as the $I_n$-maximal ones,
since we have assumed that no relations are added between the elements of $I_n$ in $I$. Therefore $a \le b$ in $\otimes^{\DS}_{i\in I}(Q^i)$.
$\qed_{\ref{lem:trecared}}$
\end{proof}

The conclusion follows by substituting $P_n^i$ for $Q_i$ in Lemma \ref{lem:trecared}.
$\qed_{\ref{generalisedDressShiffels}}$
\end{proof-}

\section{Conclusion}\label{sec: conclusion} In this paper we have been motivated by a tempting conjecture, formulated in \cite{bqoconjecture}, which to some extent claims that all wqo are ``morally'' speaking bqo. Namely, that even though the Rado example is a well known example of a wqo which is not bqo, other examples are hard to build and do not tend to appear naturally. The conjecture states this in a more precise way, claiming that every wqo is a countable union of bqos. The Rado example obviously satisfies the conjecture. Although we have not solved the conjecture, we made progress in confirming that if a counterexample is to be found, then it is not to be done by operations on wqo that have appeared in the literature. 

In searching for examples of operations on wqos, we have rediscovered a little known, but it seems to us interesting, operation which we named the Dress-Schiffels product.
Indeed, this product is defined by what one could name a finite support product of wqo posets. The operation was introduced in technical report \cite{DressSchiffels}, but apart from the paper \cite{MR1440456}, it has not appeared ever again. The known results about it were reproved in an unpublished book manuscript \cite{Mauricebook}. In Section \ref{sec:DS}  we have taken it from there and developed a number of properties of this product. 

We also considered the ordinary products of posets and have developed a set of results connecting the properties of posets of the form $P^\alpha$ and $P^{<\alpha}$, Section
\ref{sec: connections}.
Interestingly, our proofs have to pass through an intermediary of the order $I(P^{<\alpha})$ of initial segments. To us that these results seem of independent interest and connect to the well-researched area of $\alpha$-wqos. Section \ref{sec:Hausdorff} is devoted only to them, giving a Hausdorff-style classification result which shows how they can be used to construct in an orderly way a larger class of posets, called $\alpha$-FAC.

Our methods often use the combinatorics of barriers developed by Nash-Williams in his work on bqos. This material is sometimes considered too technical and difficult to remember. In Section \ref{reviewbqo} we have given a presentation which presents this material as an easily memorable extension of the well known Ramsey theorem. We hope that this presentation will make the methods more usable to a wider audience.


\bibliographystyle{splncs04}
\bibliography{../../../bibliomaster}

\end{document}